\documentclass[9pt]{article}

\usepackage{amsmath}
\usepackage{amsfonts}
\usepackage{amscd}
\usepackage{amsthm}
\usepackage{multirow}

\newtheorem{theorem}{Theorem}
\newtheorem{cor}{Corollary}

\begin{document}

\title{The Splitting of Primes in Division Fields of
 Elliptic Curves  }

\author{W.Duke
\thanks{Supported by NSF grant DMS-98-01642, the Clay
Mathematics Institute and the American Institute of Mathematics}
and \'A. T\'oth \thanks{Supported by a Rackham grant}
}

\date{}
\maketitle

\begin{center}
{\it Dedicated to the memory of Petr \~{C}i\v{z}ek}
\end{center}

\begin{abstract}
In this paper we will give a global
description of the Frobenius for the division fields of an elliptic
curve $E$ which is strictly analogous to the cyclotomic case. This is
then applied to determine the splitting of a prime $p$ in
subfields of such a division field.
Such fields include a large class of non-solvable quintic
extensions and our application provides an arithmetic
counterpart to Klein's ``solution'' of quintic equations using elliptic
functions.
A central role is played by the discriminant of the ring of
endomorphisms of the elliptic curve reduced modulo $p$.
\end{abstract}

                      %%%%%%%%%%%%%%%%%%%%%%

%%%%%%%%%%%%%%%%%%%%%%     SECTION: INTRO

                      %%%%%%%%%%%%%%%%%%%%%%

\section*{Introduction}

Given a Galois extension $L/K$ of number fields with Galois group $G$,
a fundamental problem is to describe the (unramified) primes
${\mathfrak p}$ of $K$ whose Frobenius automorphisms lie in a given
conjugacy class $C$ of $G$. In particular, all such primes have the
same splitting type in a sub-extension of $L/K$. In general, all that
is known is that the primes have density $|C|/|G|$ in the set of all
primes, this being the Chebotarev theorem.

For $L/K$ an abelian extension Artin reciprocity describes such primes
by means of their residues in generalized ideal classes of $K$.  In
the special case that $L$ is obtained explicitly by adjoining to $K$
the $q$-th division points of the unit circle we have that $G \subset
GL(1,q)=(\mathbb{Z}/q\mathbb{Z})^*$ and the Frobenius of ${\mathfrak
p}$ is determined by the norm $N({\mathfrak p})$ modulo $q$.  If
$K=\mathbb{Q}$ (cyclotomic fields) we have that $G = GL(1,q)$ and any
abelian extension of $\mathbb{Q}$ occurs as a subfield of such an $L$
for a suitable $q$ (Kronecker-Weber). Here the Chebotarev theorem
reduces to the prime number theorem in arithmetic progressions.

In a similar manner an elliptic curve $E$ over $K$ gives rise to its
$q$-th division field $L_q$ by adjoining to $K$ all the coordinates of
the $q$-torsion points.  Now $L_q$ is a (generally non-abelian) Galois
extension of $K$ with Galois group $G$, a subgroup of
$GL(2,\mathbb{Z}/q\mathbb{Z})$ (see \cite{Serre}).  In this paper we will give a global
description of the Frobenius for the division fields of an elliptic
curve $E$ which is strictly analogous to the cyclotomic case. This is
then applied to determine the splitting of primes in fields contained
in $L_q$ or, as we shall say, uniformized by $E$. As observed by
Klein (see \cite{Klein}), such fields include a large class of non-solvable quintic
extensions.  Our aim in this application is to provide an arithmetic
counterpart to Klein's ``solution'' of quintic equations using elliptic
functions.

By using CM curves we may uniformize all abelian extensions of
imaginary quadratic fields.  A classical application here is the
result of Gauss that
$$
x^3-2
$$ factors completely modulo a prime $p>3$ if and only if
$$p=x^2+27y^2$$
for integers $x$ and $y$ (see \cite{Cox}).
One way to derive this is to determine the Frobenius class of $p$ in the field obtained by
 adjoining to $\mathbb{Q}$ the x-coordinates of the 3-division points of the elliptic curve given by
$$y^2=x^3-15x+22,$$
which has CM by the quadratic order of discriminant -12.

Analogous results for non-solvable quintics require non-CM curves.
Consider the quintic $$f(x)=x^5+90x^3+3645x-6480$$ which has
discriminant $(2)^{12}(3)^{16}(5)^5(7)^6$. Its splitting field has
Galois group $S_5$ over $\mathbb{Q}$.  It follows from the results of
this paper that $f(x)$ factors completely modulo $p >7$ if and only if
$$p=x^2-25 \Delta_p y^2$$ where $\Delta_p$ is the discriminant of the ring of
endomorphisms of the elliptic curve
$$ y^2=x(x-1)(x-3)$$ reduced mod $p$.  The first two such primes are
1259 and 1951 for which $\Delta_{1259}=-31$ and $\Delta_{1951}=-51$
and where
$$1259=(22)^2+25\cdot31\cdot1^2$$
and
$$1951=(26)^2+25\cdot51\cdot1^2.$$
As may be checked,
$$f(x) \equiv (x+734)(x+322)(x+26)(x+851)(x+585)\bmod 1259$$
and
$$f(x) \equiv (x+1029)(x+1222)(x+839)(x+1771)(x+992) \bmod 1951.$$
In the non-CM case $\Delta_p$ is not determined by
arithmetic progressions in $p$. A goal of this paper is to complement that of Shimura \cite{Shimura}
by pointing out the role of $\Delta_p$ in such questions.

\bigskip
\noindent
{\sc Acknowledgement:} We would like to thank N. Katz for his helpful comments.

\section*{Outline of results.}\label{}

Given an elliptic curve $E$ defined over a number field $K$ and a
prime ideal ${\mathfrak p}$ in $\cal{O}_{K}$ of good reduction for $E$
we shall define an integral matrix $[{\mathfrak p}]$ of determinant
$N({\mathfrak p})$ whose reduction modulo $q$ gives the action of the
Frobenius for $L_q$, the $q$-th division field of $E$.  Let
$a_{\mathfrak p}$ be defined as usual by
\begin{equation}\label{1}
\#E_{\mathfrak p}(k) = N({\mathfrak p}) - a_{\mathfrak p} + 1
\end{equation}
\noindent
where $E_{{\mathfrak p}}$ is the reduction of $E$ at $\mathfrak p$ and
is defined over $k$, the residue field of ${\mathfrak p}$
which satisfies $\#k=N({\mathfrak p})= p^n.$

Let $R$ be the ring of those endomorphisms of $E$ that are rational
polynomial expressions of the Frobenius endomorphism $\phi_k$. If $\phi_p$ is
multiplication by an integer $R=\mathbb{Z}$ and we define
$\Delta_{\mathfrak p}=1$ and $b_{\mathfrak p}=0$. Otherwise the ring
$R$ is the centralizer of the Frobenius endomorphism in the
endomorphism ring of $E_{\mathfrak p}$ over $k$ and is an imaginary
quadratic order whose discriminant we denote by $\Delta_{\mathfrak
p}$. Then we shall see that $p$ does not divide the conductor $m$ of
$\Delta_{\mathfrak p}$ and that there is a unique non-negative integer
$b_{\mathfrak p}$ so that
\begin{equation}\label{2}
4N({\mathfrak p}) \;=\; a_{\mathfrak p}^2 \;-\;\Delta _{\mathfrak p}^{ }\;b_{\mathfrak p}^2.
\end{equation}
We associate to ${\mathfrak p}$ the following integral matrix of determinant $N({\mathfrak p})$:

\begin{equation}\label{3}
[{\mathfrak p}] = \begin{bmatrix} (a_{\mathfrak p} + b_{\mathfrak p} \delta_{\mathfrak p})/2 &
b_{\mathfrak p}
\\
b_{\mathfrak p}(\Delta_{\mathfrak p}- \delta_{\mathfrak p})/4 & (a_{\mathfrak p} - b_{\mathfrak p}
\delta_{\mathfrak p} )/2
\end{bmatrix}
\end{equation}
\noindent
where for a discriminant $\Delta$ we have $\delta = 0,1$ according to
whether $\Delta \equiv 0,1\bmod 4$. We shall show that $[{\mathfrak
p}]$ gives a global representation of the Frobenius class over
${\mathfrak p}$ for each $q$-th division field of $E$ by reducing it
modulo $q$, provided $p$ is prime to $q$.

\medskip
\noindent

\begin{theorem}\label{Thm1}
Let $E$ be an elliptic curve defined over a number field $K$ and $q>1$ an
integer.  Let $L_q$ be the $q$-th division field of $E$ with Galois group $G$ over $K$.
Let ${\mathfrak p}$ be a prime of good reduction for $E$ with $N(\mathfrak p)$
prime to $q$. Then ${\mathfrak p}$
is unramified in $L_q$  and the
integral matrix $[{\mathfrak p}]$ defined in (\ref{3}), when reduced modulo
$q$, represents the class of the Frobenius of ${\mathfrak p}$ in $G$.
\end{theorem}

\noindent
The proof we give of this uses the theory of canonical lifts of endomorphisms
due originally to Deuring.

In analogy with the cyclotomic case, we have associated to each curve
a sequence of prime power matrices, defined in terms of arithmetic
data from the reduced elliptic curve which give the Frobenius in all
of the division fields.  Let $C$ be a conjugacy class of $G$ and let
$\pi_E(X;q,C)$ be the number of primes ${\mathfrak p}$ of good
reduction with $N({\mathfrak p})\leq X$ such that $[{\mathfrak p}]
\equiv C_0~\bmod q$ for some $C_0 \in C.$ By the Chebotarev theorem \cite{Cheb}
we
derive the following strict analogue of the prime number theorem in
progressions for the sequence $[{\mathfrak p}]:$
$$\pi_E(X;q,C) \sim \frac{|C|}{|G|} \pi_K(X)$$ as $X \rightarrow
\infty$, where $\pi_K(X)$ counts all primes of $K$ with $N({\mathfrak p})\leq X.$

Of more interest for us here is the fact that the splitting type of
 ${\mathfrak p}$ in any field between $K$ and the $q$-th division
 field $L_q$ is determined by $[{\mathfrak p}] \bmod q$. For example
 we get immediately a criterion for complete splitting in the full
 division field in terms of the invariants $a_{\mathfrak p}$ and
 $b_{\mathfrak p}$ modulo $q$, provided $q$ is odd.

\begin{cor}\label{cor1}
Let $E$ be an elliptic curve defined over a number field $K$ and $q>1$ an
odd integer.
Then ${\mathfrak p}$ a prime of good reduction for $E$ with $N(\mathfrak p)$
prime to $q$
splits completely in
 $L_q$ if and only if  $a_{\mathfrak p} \equiv 2
\bmod q $ and $b_{\mathfrak p} \equiv 0 \bmod q $.
\end{cor}
\bigskip

For a discriminant $\Delta$ let
$$Q_{\Delta}(x,y)= x^2 + \delta x y -((\Delta-\delta)/4)y^2$$
be the principal form where $\delta = 0,1$ according to whether $\Delta \equiv 0,1\bmod
4$. For ${\mathfrak p}$ a prime of good reduction for $E$
we get a representation
\begin{equation}\label{4}
N({\mathfrak p})= Q_{\Delta_{\mathfrak p}}(x,y)
\end{equation}
with integral $x,y$ upon using the change of variables
\begin{equation}\label{5}
x=(a_{\mathfrak p}-b_{\mathfrak p} \delta_{\mathfrak p})/2 \,\,\,\,\,\,\,\,\,\,\,\,
y=b_{\mathfrak p}
\end{equation}
in (2). This representation is primitive if $\mathfrak p$ is ordinary.
Let $L_q^+$ be the extension of $K$ obtained by adjoining only the
Weber functions of the $q$-th division points, that is the
$x$-coordinates unless $j(E)=0$ or $j(E)=1728$, in which case we must
first cube or square the coordinates, respectively.  By Theorem 1 we
may determine which sufficiently large ordinary primes split
completely in $L_q^+$ from any such primitive representation.

\begin{cor}\label{cor3}
Let $E$ be an elliptic curve defined over a number field $K$ as above
and $q \geq 1$ an integer.  Then there is a constant $C_0$ depending
only on $E$ and $q$ so that for every ordinary prime ${\mathfrak p}$
of $K$ with $N({\mathfrak p})>C_0$ we have that $\mathfrak p$ splits
completely in $L_q^+$ if and only if $x \equiv \pm 1 \bmod q$ and $y
\equiv 0 \bmod q$ in any primitive representation
$$N({\mathfrak p})= Q_{\Delta_{\mathfrak p}}(x,y).$$
\end{cor}

\bigskip
If $E$ has CM by the ring of integers in an imaginary quadratic field
of discriminant $\Delta$ then the splitting completely condition in
$L_q^+$ becomes simply
$$N({\mathfrak p})= Q_{\Delta}(x,y)$$ with integers $x \equiv \pm 1
\bmod q$ and $y \equiv 0 \bmod q$.  Actually, suppose we take for $E$
the elliptic curve with lattice given by the ring of integers of an
imaginary quadratic field $F$ of discriminant $\Delta$ and take
$K=F(j(E))$, the Hilbert class field of $F$. It follows from Corollary
2 that a sufficiently large rational prime $p$ splits in $L_q^+$ iff
$p=Q_{\Delta}(x,y)$ with integers $x \equiv \pm 1 \bmod q$ and $y
\equiv 0 \bmod q$.  This is a well known result of CM theory.

Another simple consequence in the CM case, this time of Corollary 1,
is that the conditions
$$\#E_{\mathfrak p}(k) \equiv 0 \bmod q^2  \,\, and \,\,
N({\mathfrak p}) \equiv 1 \bmod q,$$
which are clearly necessary  for ${\mathfrak p}$ of good reduction to
split completely in $L_q$,
are also sufficient, at least when $q$ is odd.

Our main application is to describe the primes which split completely
in certain non-solvable quintic extensions  $M/K$.
Suppose $M$ is given by adjoining to $K$ a solution of
a principal quintic over $K$:
$$f(x)=x^5+ax^2+bx+c=0$$
and that the discriminant of $f$ is $D$.
Suppose further that the Galois group of
the normal closure $L$ of $M$ is $S_5$
and that $\sqrt{5D} \in K.$

\begin{theorem}\label{splitting-quintic}
Let $M/K$ be a non-solvable quintic extension as above.
There exists an elliptic curve $E$ defined over $K$  so that a prime ${\mathfrak p}$ of $K$
which has good reduction for $E$
and is prime to 5 splits completely in $M$ if and only if

$$
b_{\mathfrak p}
\equiv 0 \bmod 5
$$
where $b_{\mathfrak p}$
is associated to the elliptic curve $E$.
\end{theorem}

In general we have the following
determination of the splitting type of $\mathfrak p$:

\bigskip
\begin{tabular}{|c||r|r|c|}
\hline
$\textrm{Splitting type of} \,\,\,\,\mathfrak p \,\,\textrm{in} \,\,\ M$ &{$\displaystyle{\left(\frac{a_{\mathfrak p}^2-4{N(\mathfrak p})}{5}\right)}$} &
{ $\displaystyle{\left(\frac{N(\mathfrak p)}{5}\right)}$}& \\
\hline
\hline
$(1)(2)^2$ & $1\ \ \ \ \ $ & $1\ \ \ $ & \\
\hline
$(1)(4)$ & $1\ \ \ \ \ $ & $-1\ \ \ $ & \\
\hline
$(1)^2(3)$ &$ -1\ \ \ \ \ $ & $1\ \ \ $ & \\
\hline
$(1)^3(2)$ &$ -1\ \ \ \ \ $ & $-1\ \ \ $ & if $ 5\vert a_{\mathfrak p}$\\
\hline
$(2)(3)$ & $-1\ \ \ \ \ $ & $-1\ \ \ $ &  if $ 5\not\vert a_{\mathfrak p}$\\
\hline
$(5)$ & $0\ \ \ \ \ $ &  \multicolumn{2}{|c|}{ if $ 5\not\vert
b_{\mathfrak p}\ \ \ $} \\
\hline
$(1)^5$ & $0\ \ \ \ \ $ & \multicolumn{2}{|c|}{ if $ 5\vert
b_{\mathfrak p}\ \ \ $}  \\
\hline

\end{tabular}

\bigskip
\noindent
Concerning the determination of $E$ from $f$, it is enough to
find the $j$-invariant of $E$. Explicit computations
are provided below.
We remark that it is also possible to formulate a similar result for $A_5$
extensions of $K$ under otherwise identical assumptions. Furthermore,
by allowing the elliptic curve to be defined over a quadratic or a biquadratic
extension of $K$ one may uniformize all non-solvable quintic extensions.

It is also possible to explicitly uniformize certain degree 7 extensions
whose normal closure have Galois group simple of order 168
by using the seventh division fields of elliptic curves
(see \cite{Radford} and the references cited there.)
By Theorem 1 one may similarly characterize the primes with a given
splitting type in such extensions.

                      %%%%%%%%%%%%%%%%%%%%%%

%%%%%%%%%%%%%%%%%%%%%%     NEW SECTION

                      %%%%%%%%%%%%%%%%%%%%%%

\section*{A global representation of the Frobenius}

In this section we will prove Theorem 1 and its corollaries using an
approach that compares the action of the Frobenius on the prime-to $p$
division points with the action of the matrix (\ref{3}) on
$\mathbb{Z}^2$.

Let $E$ be an elliptic curve defined over a number field $K$. Let
$\mathfrak p$ be a prime ideal in $\cal{O}_{K}$ with residue field $k$
$E_{\mathfrak p}$, the reduction of $E$ mod $\mathfrak p$ (it is assumed
that $E$ has good reduction at ${\mathfrak p}$).  That ${\mathfrak p}$
is unramified in the field $L_q$ is well known, see
e.g. \cite{Silverman} VII.\S4. Also note that there is nothing to
prove when $\phi_{\mathfrak p} \in \mathbb{Z}$, so we will assume
throughout that this is not the case.  The idea of the proof is that
modulo ${\mathfrak p}$ the curve $E$ can be replaced by a curve
$\tilde{E}$ with complex multiplication so that the following diagram
commutes:

\noindent
\begin{equation}\label{diagram}
\begin{CD}
E[q] @>{red}>> E_{\mathfrak p}[q] @<{red}<< \tilde{E}[q]  \\
@V{F_{\mathfrak P}}VV   @V{\phi_{\mathfrak p}}VV
@V{\tilde{\phi}_{\mathfrak p}}VV \\
E[q] @>{red}>> E_{\mathfrak p}[q] @<{red}<< \tilde{E}[q]   \\
\end{CD}
\end{equation}

\noindent
where as usual $[q]$ stands for the $q$-division points on the curves
in the algebraic closures of the appropriate fields.

We now explain this diagram in detail. To simplify matters we fix a
Weierstrass equation for $E$ as in \cite{Silverman} III.\S1. Let
$\overline{K}, \overline{k}$ be the algebraic closures of $K, k$. To
specify the horizontal maps $red$ we choose an embedding of
$\overline{K}$ into the algebraic closure $\overline{K_{\mathfrak p}}$
of $K_{\mathfrak p}$, the completion of $K$ at the valuation arising
from ${\mathfrak p}$.  We call the subgroup of torsion points whose
orders are relatively prime to $p$ the $p'$-torsion. Then the
$p'$-torsion points on $E(\overline{K})$ is mapped into the
$p'$-torsion of $E(\overline{K_{\mathfrak p}})$ and this being defined
over an unramified extension, reduction modulo a prime ${\mathfrak P}$
above ${\mathfrak p}$ maps this latter group into the $p'$-torsion of
$E(\overline{k})$.  Both of these maps are isomorphisms on $p'$
torsion. This is the map $red$ for reduction, though as explained
above it depends on many choices. Note that after these choices are
made there is a unique element $F_{\mathfrak p} \in Gal(K_{\mathfrak
p}^{unram}/K_{\mathfrak p})$ that satisfies $F_{\mathfrak p}(t) \equiv
t^{\#k} \bmod {\mathfrak P}$, for all $t\in K_{\mathfrak p}^{unram}$.

We are interested in the action of the Frobenius automorphism
$\phi_{\mathfrak p} \in Gal(\overline{k}/k)$ on the
$\overline{k}$-valued points. In terms of the Weierstrass equation for
$E$, this action on the coordinates is simply $(x,y) \mapsto
(x^{\#k},y^{\#k})$. By abuse of notation we also denote this action
and the restriction of it to the $q$-division points by
$\phi_{\mathfrak p}$.

\medskip
Now the commutativity of the left half of the diagram is merely a
restatement of the choices made above.

By Deuring's lifting theorem (\cite{Deuring}, \cite{Lang} p.184),
there exists an elliptic curve $\tilde{E}$ defined over ${K_{\mathfrak
p}}$ and an endomorphism $\tilde{\phi}_{\mathfrak p}$ of $\tilde{E}$
so that $\tilde{E}$ reduces to $E_{\mathfrak p}$ modulo ${\mathfrak
p}{\cal O}_{\mathfrak p}$ and that $\tilde{\phi}_{\mathfrak p} \in
End(\tilde{E}) $ reduces to $\phi_{\mathfrak p} \in End(E_{\mathfrak
p})$. If $E$ is super-singular $\tilde{\phi}_{\mathfrak p}$ will be
defined over a ramified extension. Reduction still makes sense since
$\tilde{\phi}_{\mathfrak p}$ is an endomorphism and not a Galois
automorphism.

\medskip
This shows the commutativity of the right half of diagram (\ref{diagram}).

\bigskip
To prove our theorem we need to determine the endomorphism ring $S$
of $\tilde{E}$.  Recall that the ring $R_{\mathfrak p}$ defined in the
introduction is the centralizer of $\phi_{\mathfrak p}$ in the
endomorphism ring of $E_{\mathfrak p}$ and is a quadratic order. We
claim that $S$ is isomorphic to $R$. Since $R \subset S$, Deuring's
reduction theorem implies equality if we can show that the conductor
of $R$ is prime to $N({\mathfrak p})$, a fact that is trivial in the
ordinary case and follows from \cite{Wat} in the super-singular case.

\medskip
Let $\Delta_{\mathfrak p}$ be the discriminant of $R_{\mathfrak
p}$.  By choosing a complex square root of $\Delta_{\mathfrak p}$ we
identify $R_{\mathfrak p}$ with a lattice in $\mathbb{C}$. After this
identification $\phi_{\mathfrak p}$ corresponds to some complex number
$\phi=(a_{\mathfrak p}+b_{\mathfrak p} \sqrt{\Delta_{\mathfrak
p}})/2$. Clearly the lattice $R$ is preserved by multiplication by
$\phi$ and leads to the integral matrix (3), where we may choose
$b_{\mathfrak p} \geq 0$.  Instead of $R$ one could in fact use any
lattice whose endomorphism ring is $R$.

To finish the proof of Theorem 1, consider an embedding $\alpha$ of
the algebraic closure of $K_{\mathfrak p}$ into $\mathbb{C}$. It
allows us to view $\tilde{E}$ as an elliptic curve over the complex
numbers, that we denote $E_{\alpha}$. Since $E_{\alpha}$ has complex
multiplication by $R$ and $Gal(\mathbb{C}/\mathbb{Q})$ acts
transitively on the set of elliptic curves with $R$ as its
endomorphism ring, we may and will assume the $j(E_{\alpha})=j(R)$.

By choosing a non-trivial holomorphic differential $\omega$ on
$\tilde{E}_\alpha$ appropriately the lattice of periods $\{
\int_\gamma \omega : \gamma \in H_1(\tilde{E}_\alpha, \mathbb{Z}
\}=R$.  Then the period mapping $ \Pi:\tilde{E}_\alpha \rightarrow
\mathbb{C}/R$ is a biholomorphic isomorphism of complex analytic
manifolds. The action of $\tilde{\phi}_{\mathfrak p}$ on $\tilde{E}$
defines an endomorphism of $\tilde{E}_\alpha$ and gives rise to a map
$\phi_*$ on $R$. Since the Frobenius automorphism $\phi_{\mathfrak p}$
satisfies a quadratic equation

\begin{equation}\label{frobenius}
\phi_{\mathfrak p}^2 - a_{\mathfrak p}\;\phi_{\mathfrak p} + N({\mathfrak p}) = 0.
\end{equation}

\noindent
$\phi_*$ can be identified with multiplication by one of the complex roots of
this equation i.e. multiplication by $\phi: R \rightarrow R$ (viewed as complex
numbers). Getting back to the $q$-division points we can again summarize the
situation in the following diagram:

\begin{equation}\label{diagram2}
\begin{CD}
\tilde{E}_{\mathfrak p}[q] @>{\alpha}>>
\tilde{E}_\alpha[q] @>q \times {\Pi}>> R /qR \\
 @V{\tilde{\phi}_{\mathfrak p}}VV
@V{\phi_\alpha}VV  @VV{\phi_*}V \\
 \tilde{E}_{\mathfrak p}[q] @>{\alpha}>>
\tilde{E}_\alpha[q] @>q \times{\Pi}>> R/qR  \\
\end{CD}
\end{equation}

\noindent
where $q \times{\Pi}$ is the period map followed by multiplication by
$q$. This proves Theorem 1.

\bigskip
\noindent
{\it Remark:} If $E$ is replaced by an Abelian variety $V$ then the
unramifiedness of $\mathfrak{p}$ still holds \cite{Shimura-Taniyama}
and the left square of diagram \ref{diagram} makes sense. If in addition
$V$ has ordinary reduction at ${\mathfrak p}$ then the right square in
diagram \ref{diagram} generalizes as shown by Deligne \cite{Del} (and
therefore the whole proof works). However the general case leads to
substantial difficulties \cite{Oort}.

\bigskip
\noindent
Corollary 1 is an immediate consequence of Theorem 1.

\bigskip
\noindent

We now prove Corollary 2.  Let $E$ be an elliptic curve defined over a
number field $K$ as above and $q \geq 1$ an integer.  Let ${\mathfrak
p}$ be a prime of ordinary reduction for $E$.  Given a primitive
representation
$$p^n= Q_{\Delta_{\mathfrak p}}(x,y)$$
we know that $x$ and $y$
are uniquely determined up to (proper or improper)
automorphs of $Q_{\Delta_{\mathfrak p}}.$
If $-\Delta_{\mathfrak p}>4$ and $x \equiv \pm 1 \bmod q$ and $y \equiv 0 \bmod q$
then
it follows that
\begin{equation}\label{7}
[{\mathfrak p}] \equiv \begin{bmatrix} x+ \delta y  & y
\\ y (\Delta_{\mathfrak p}- \delta_{\mathfrak p})/4 & x
\end{bmatrix} \bmod q
\end{equation}
and hence that $\mathfrak p$ splits completely in $L_q^+$.  If
$j=j(E)$ is not 0 or 1728 then for $\mathfrak p$ with $N(\mathfrak p)$
sufficiently large we have that $-\Delta_{\mathfrak p}>4$. To see this
write $j=\alpha/\beta$ for $\alpha,\beta \in \cal O_K$. We know
that $j \equiv j(R_{\mathfrak p}) \bmod \mathfrak p$.  If
$j(R_{\mathfrak p})=0$ or 1728 then assuming that $j - j(R_{\mathfrak
p}) \neq 0$ we have
$$N(\mathfrak p) \leq \max(|N(\alpha)|,|N(\alpha-1728\beta)|).$$
In  case  $j=0$ or $j=1728$ the altered definition of $L_q^+$
leads again to the result.

\bigskip
\noindent

Finally we prove the consequence of Corollary 1  mentioned below Corollary 2 above that, in the CM case,
a prime of good reduction ${\mathfrak p}$ splits completely in $L_q$ if
$a_{\mathfrak p} \equiv N({\mathfrak p})+1 \bmod q^2 $ and
$N({\mathfrak p}) \equiv 1 \bmod q$, provided $q$ is odd.
Since these conditions immediately imply that $a_{\mathfrak p} \equiv 2 \bmod{q},$
by Corollary 1 we only must show that $q \mid b_{\mathfrak p}$.
By our assumption
$$
a_p^2 \equiv (N({\mathfrak p})-1)^2 +4N({\mathfrak p})
\equiv 4N({\mathfrak p}) \bmod{q^2}
$$
we get,  using
$$4N({\mathfrak p}) \;=\; a_{\mathfrak p}^2 \;-\;\Delta_{\mathfrak p} \;b_{\mathfrak p}^2,$$
that
$$q^2 \mid \Delta_{\mathfrak p} b_{\mathfrak p}^2.$$
For a CM curve with fundamental $\Delta$ the only possible prime
dividing the square part of $\Delta_{\mathfrak p}$ is 2.  In fact,
$\Delta_{\mathfrak p}= \Delta$ for ordinary $\mathfrak p$ and for
super-singular $\mathfrak p$ $\Delta_{\mathfrak p}=-p$ or $-4p$ where
$N(\mathfrak p)=p^n.$ Since $q$ is odd this implies that $q \mid
b_{\mathfrak p}$.
\bigskip
\noindent

                      %%%%%%%%%%%%%%%%%%%%%%

%%%%%%%%%%%%%%%%%%%%%%     NEW SECTION

                      %%%%%%%%%%%%%%%%%%%%%%

\section*{Quintics}

In the section we prove Theorem 2 and justify the general splitting criteria
given after it as well as the example given in the introduction.
Let $M$ be given by adjoining to $K$ a root of
a principal quintic
$$f(x)=x^5+ax^2+bx+c=0$$
defined over $K$.
If the discriminant of $f$ is $5$ times a square
then, by means of a Tschirnhausen transformation
(\cite{Dickson} p.218.)
we may assume that $M$ is determined by
a Brioschi quintic
$$
f_t(x)=x^5-10tx^3 + 45t^2x-t^2
$$
for some $t \in K$ with $t \neq 0, \frac{1}{1728}$ .
It was shown by Kiepert \cite{Kiepert} already in 1879 (see \cite{King} for an exposition)
that $M$ is contained in $L_5^+$ for any elliptic
curve $E$ over $K$ with $j$-invariant $1728-t^{-1}$.
Recall that $L_5^+$ is in this case obtained by adjoining to $K$ the $x$-coordinates
of the 5 division points.
One may take for instance the curve $E_t$ given by
$$
E_t: y^2+xy = x^3+36tx+t.
$$

If the splitting field of $f$ over $K$
is an $S_5$ extension then it must be the fixed field of the
subgroup of scalars of $G$ since $PGL_2(\mathbb{F}_5) \simeq S_5$.
Theorem 2 now follows easily from Theorem 1.

A calculation of conjugacy classes based on the identification of $S_5$
with $PGL_2(\mathbb{F}_5)$ leads to the determination of the splitting
type of a prime $\mathfrak p$ of good reduction for $E_t$ which is prime to 5.
Recall that $A \in GL_2(\mathbb{F}_5)$ is called regular if it has
different eigenvalues.
Clearly $A$ is regular if the discriminant of the characteristic
equation $tr(A)^2-4\det(A)$ is non-zero. Given such  $A$
its conjugacy class is determined by its trace and determinant.
It is clear that the values of the following Legendre symbols
$$
\sigma=\left(\frac{\det(A)}{5}\right)\ \ \ \text{ and } \ \ \
\rho=\left(\frac{tr(A)^2-4\det(A)}{5}\right)
$$
\noindent
are determined by the conjugacy class of $A$ in $PGL_2(\mathbb{F}_5)$. Now in
case the characteristic polynomial of $A$ splits, that is $\rho=1$, the matrix $A$ is
conjugate to a diagonal matrix in $GL_2(\mathbb{F}_5)$ and so the value of
$\sigma$ already determines the cycle type of
such matrices. When $\rho=-1$ one
must take into account whether $tr(A)\equiv 0 \text{ or } \not \equiv 0
\bmod 5$. For $A$ non-regular
$tr(A)^2 -4 \det(A)=0$ and one needs to know if $A$ is semi-simple or
unipotent. This information cannot be extracted from the trace and
determinant alone, but it is determined by the value of $b_{\mathfrak p}$.
All that remains to be done is to identify each conjugacy classes with its
cycle type.

The example in the introduction is obtained by taking $K=\mathbb{Q}$
and $t=\frac{-3^2}{2^85^2} $. Here we observe that since $E$ has four
2-torsion points over $\mathbb{Q}$ both $a_p$ and $b_p$ will be even
for $p$ with good reduction. Thus the representation
$$4p \;=\; a_{p}^2 \;-\;\Delta _{p}\;b_{p}^2$$
yields
$$p \;=\; x^2 \;-\;\Delta _{p}\;y^2$$
and the condition for splitting completely
is that $y \equiv 0 \bmod 5,$ since $x$ and $y$ are determined uniquely
up to sign.

\section*{Some computational issues}

In this section we discuss some of the computational issues
which arise when considering examples.

First, given a principal quintic (slightly modified from above)
\begin{equation}\label{principal}
f(x)=x^5+5ax^2+5bx+c=0
\end{equation}
\noindent
defined over $K$ with discriminant $D$ such that $\sqrt{5D} \in K$
we must determine $t$ so that
the Brioschi quintic
\begin{equation}\label{brioschi}
f_t(x)=x^5-10tx^3 + 45t^2x-t^2.
\end{equation}
\noindent
determines the same extension.
This is done using a Tschirnhausen transformation and is described in detail
in \cite{King}, p 103. (see also \cite{Dickson} p.218.)
Here we will simply record the result in the case $a \neq 0.$

One determines $t, \lambda$ and $\mu$ in the map
\begin{equation}\label{tschirn}
x\mapsto \frac{\lambda+\mu x}{(x^2/t) -3}
\end{equation}
\noindent
in order to transform the general principal quintic (\ref{principal}) to the Brioschi quintic (\ref{brioschi}).

An analysis using invariant polynomials for the icosahedral group acting on the Riemann sphere
leads eventually to the quadratic equation for $\lambda$
given by
$$ (a^4+abc-b^3)\lambda^2-(11a^3-ac^2+2b^2c)\lambda+64a^2b^2-27a^3c-bc^2=0.$$
The discriminant of this quadratic is
$$5^{-5}a^2D$$
and so $\lambda \in K$. Choose either solution and
let
$$j=\frac{(a\lambda^2-3b\lambda-3c)^3}{a^2(\lambda ac-\lambda b^2 -bc)}.$$
Then, provided $j \neq 0, 1728$ we may take
$$t=1/(1728-j)$$
 in (\ref{brioschi})
and choose for the elliptic curve any curve with this $j$ invariant,
say
$$
E_t: y^2+xy = x^3+36tx+t
$$
as above.
Also, one may determine
$\mu$ in (\ref{tschirn}) to be given by
$$ \mu=\frac{j a^2-8{\lambda}^3a-72{\lambda}^2b-72\lambda c}{{\lambda}^2a+\lambda b +c}.$$
Note that the discriminant of $f_t$ is $$D_t=5^5t^8(1728t-1)^2$$ while that of
$E_t$ is $$-t(1728t-1)^2.$$

Another issue is to compute the invariants $\Delta_p$ and $b_p$ in the
rational case. This problem is mentioned in the Woods Hole seminar on
Elliptic curves and formal groups by Lubin-Serre-Tate.
An important study of $\Delta_p$ was made by Schoof in \cite{Schoof}.

The most straightforward way to determine $b_p$ and to find the
order $R$ that appears in Deuring's theorem is to
check all the possible singular invariants until we find one that is
congruent to the given $j$-value modulo $p$.  (Note that the
discriminant of $R$ must divide $a_p^2-4p$.) We assume that our input
is an elliptic curve $E$, given in Weierstrass equation, and $p$ is a
prime number that does not divide the discriminant of $E$. After
computing $a_p$ (a serious task by itself) we find $\Delta_p$ for an
ordinary curve as follows; we first compute the square-free part $D$
of $a_p^2-4p$ and then create a vector whose values are all possible
discriminants
$$\Delta=b^2 D|(a_p^2-4p).$$

\smallskip
\noindent
For a possible conductor $\Delta$, we find the class group ${\cal
C}(\Delta)$ of the proper ideal classes (using quadratic forms) and
compute the integer
\[
X_{\Delta}=\displaystyle{\prod_{\Lambda \in {\cal C}(\Delta)}
\left(j(E)-j(\Lambda)\right)}.
\]

\noindent
Note that the canonical lift $\tilde{E}$ is distinguished by the fact
that its endomorphism ring is $R$ and that $j(E) \equiv j(\tilde{E})
\bmod{P}$ for some prime $P$ dividing $p$. Therefore for any
complex embedding
$\alpha: \mathbb{Q}_p \rightarrow \mathbb{C},$
$$\alpha(j(\tilde{E})) \in \{j(\mathbb{C}/\Lambda): \Lambda \in {\cal
C}_{\Delta_p}\},$$
 where $\Delta_p$ is the actual discriminant of
$R$. Also note that if $\Lambda' \in {\cal C}_{\Delta'}$ for
$\Delta'\neq \Delta_p$, then the corresponding elliptic curve reduces
to a curve whose endomorphism ring has discriminant $\Delta'$ for any
place above $p$. Therefore $\Delta_p$ is uniquely characterized by the
fact that $$X_{\Delta_p} \equiv 0 \bmod p.$$

Occasionally the computation of $X_{\Delta}$ involves complex
numbers of rather large size. To make the algorithm efficient, one needs
to determine the needed precision in advance. Assume that the lattices
are given in the form $\mathbb{Z}+\mathbb{Z}\tau_i$, with $\tau_i$ in
the upper half plane. Then the number of significant digits one must
use is approximately

\[
\sum_{\tau_i}^{}
\frac{\log(j(E)) + 2 \pi Im(\tau_i)}{\log(10)}.
\]

\noindent

\newpage

\bibliographystyle{plain}

\end{document}